\documentclass[11pt]{article}
\usepackage{amsmath}
\usepackage{amsthm}
\usepackage{amssymb}
\usepackage{graphicx}

\renewcommand{\leq}{\leqslant}
\renewcommand{\geq}{\geqslant}
\renewcommand{\le}{\leqslant}
\renewcommand{\ge}{\geqslant}

\newcommand{\su}{\subseteq}

\newcommand{\midd}{:}

\newcommand{\R}{\mathbb{R}}

\newtheorem{theorem}{Theorem}[section]
\newtheorem{corollary}[theorem]{Corollary}
\newtheorem{conjecture}[theorem]{Conjecture}
\newtheorem{lemma}[theorem]{Lemma}

\newtheorem{claim}[theorem]{Claim}

\newcommand{\ignore}[1]{}
\date{}

\author{Noga Alon\thanks{Department of Mathematics,
Princeton University, Princeton, NJ 08544, USA.
Research supported in part by NSF grant DMS-2154082.
\texttt{nalon@math.princeton.edu}}
\and Rom Pinchasi\thanks{Mathematics Department, 
Technion, Haifa 32000, Israel.
Supported in part by ISF grant No. 1091/21.
\texttt{room@technion.ac.il}.}}

\title{Distinct Directions and Distinct Distances in $\mathbb{R}^d$} 

\begin{document}

\maketitle

\begin{abstract}
We show that there exists an absolute positive constant $b (\geq
\frac{1}{48})$ so that
any set of $n$ points in $\mathbb{R}^d$ that is 
$d$-dimensional determines at least $bdn$ lines 
with pairwise distinct directions. As a consequence, we prove that 
there are 
$d$-dimensional real norms $\|\cdot\|$ 
so that every set of $n>n_0(d)$ points that is 
$d$-dimensional determines at least 
$(bd-o(1))n$ distinct distances with respect to $\|\cdot \|$.
%In fact, this holds for all  $d$-norms but a meagre set.
\end{abstract}

\section{Introduction}

The celebrated Gallai-Sylvester theorem asserts that $n$ points in 
the plane that are not collinear must determine 
an \emph{ordinary line}, that is, a line passing through precisely 
two points of the set. Erd\H{o}s noticed the following simple consequence. 
Any set of $n$ points in the plane that are not collinear must 
determine at least $n$ distinct lines, and equality occurs 
only when $n-1$ of the points are collinear.

In the same spirit, Scott asked two similar questions in 1970:

\begin{enumerate}

\item  Is it true that the minimum number of distinct directions 
of lines determined by $n$ noncollinear points in $\mathbb{R}^2$ is 
$2\lfloor \frac{n}{2} \rfloor$?

\item What is the minimum number of distinct directions of 
lines determined by a $3$-dimensional set of $n$ points in $\mathbb{R}^3$?

\end{enumerate}
In 1982 Scott's first question was answered in the affirmative 
by the celebrated theorem of Ungar \cite{Un} using the technique of 
allowable sequences invented by Goodman and Pollack.

Scott's second question was answered only much later in \cite{PPS},
where it is shown that a $3$-dimensional set of $n$ points 
determines at least $2n-5$ lines with distinct directions. 
This bound is sharp when $n$ is odd. 

It is now natural to wonder what happens in higher dimensions.
The complexity of the solution of Scott's problem in $\mathbb{R}^3$, 
compared with the solution in $\mathbb{R}^2$, suggests that the
question is not likely to be easier in higher dimensions. 
Similarly, the actual 
bound in three dimensions and the sharp examples suggest that one 
cannot expect a very simple expression as a tight answer for the same 
problem in $d$ dimensions, nor a very simple construction 
achieving the minimum possible number of distinct directions.

It is however plausible to conjecture the following, which is suggested
in \cite{BMP}, page  272, motivated by earlier related results of
Jamison \cite{Ja} and of Blokhuis and Seress \cite{BS}.

\begin{conjecture}
\label{c11}
A $d$-dimensional set of $n$ points determines at 
least $(d-1)n-O(d^2)$ lines with pairwise distinct directions.
\end{conjecture}

The example of $n-d+1$ collinear points plus additional $d-1$ 
points in general position 
not on this line shows that one cannot expect the answer to be more than 
$$(d-1)(n-d+1)+{{d-1} \choose 2}+1=(d-1)n-O(d^2)$$ 
even if we wish 
to bound from below the total number of lines determined by a 
$d$-dimensional set of $n$ points, regardless of the directions 
of these lines.

It is interesting to note that the problem of finding the minimum number of distinct directions determined by $n$ points in the plane and also in higher dimension is directly related to an old conjecture of Dirac \cite{Di} in the plane. A slightly modified version of this conjecture of Dirac
asserts that for any set of $n$ points in the plane, not all on a 
single line, there is a point that lies in at least 
$n/2-O(1)$ distinct lines determined by the set. This Conjecture
is still open but several weaker versions have been established
over the years. See \cite{H17} for the best known bound and the
history of the problem. It is shown in \cite{H17} that in any two dimensional set $P$ of $n$ points one can find a point $x \in P$
that lies on at least $\frac{n}{3}$ lines determined by $P$.
Clearly, the lines determined by $P$ that are incident to $x$ have pairwise distinct directions. We will make use of this observation.

In the present note we prove a modest result supporting 
Conjecture \ref{c11} and
show that the number of distinct directions determined by 
a $d$-dimensional set of $n$ points is at least linear in $dn$.  

\begin{theorem}
\label{t12}
There exists an absolute constant $b>0$ such that for every $d \geq
2$, any set of $n$ points in $\mathbb{R}^d$  
that is not contained
in a hyperplane determines at least $bdn$ 
line segments with pairwise distinct directions.
\end{theorem}

Our proof shows that the result above holds with $b=\frac{1}{48}$. This
estimate can be easily improved, but since our method does not suffice
to obtain the best possible $b$ (which should be close to $1$ for
large $d$, if Conjecture
\ref{c11} holds), we make no attempt to optimize it.

Unlike in two and three dimensions, where the proofs of 
the sharp bounds are rather involved, our proof is shorter and is 
derived as a consequence of a deep result by Dvir, Saraf, 
and Wigderson (\cite{DSW14}), improving on an earlier similar 
result by Barak, Dvir, Wigderson, and Yehudayoff (\cite{BDWY11}).

As suggested in the solution of Scott's problem in $3$-dimensions (\cite{PPS}), it is helpful to consider pairwise non-convergent segments, rather than pairwise non-parallel segments, as the former property, unlike the latter one, is much better preserved under projective transformations. As our argument will involve central projections,
we will use a slightly generalized definition that will replace definition of two segments not being parallel. Our setting will involve a bichromatic set of points.

\begin{figure}[ht]
	\centering
	\includegraphics[width=3in]{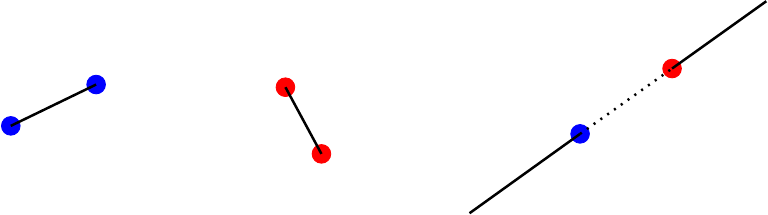}
	\caption[]
	{Generalized segments determined by red and blue points. The dotted portion on the right is not part of the generalized segment determined by the corresponding red point and blue point.}
	\label{fig:generalized}
\end{figure}

For a set of red and blue points in $\mathbb{R}^d$ every pair of points $a$ and $b$ defines a {\emph generalized segment} on the line $\ell$ through $a$ and $b$. This generalized segment is the classical closed straight line segment $[a,b]$ delimited by $a$ and $b$ if both $a$ and $b$ have the same color. However, if $a$ and $b$ have different colors, then the generalized segment that is defined by $a$ and $b$ is the set of these points on $\ell$ that do not lie between $a$ and $b$. In other words, it is equal to $\ell$ minus the open interval $(a,b)$ (see Figure \ref{fig:generalized}). The reason for this definition will become clear later.

\begin{figure}[ht]
	\centering
	\includegraphics[width=3in]{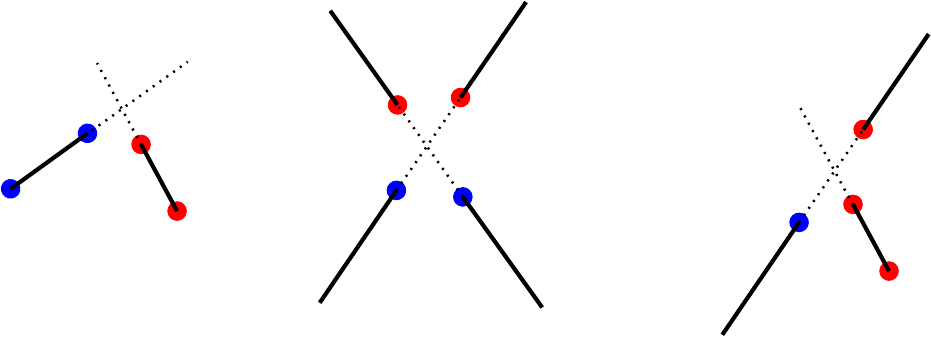}
	\caption[]
	{Pairs of (co-planar) convergent generalized segments. The dotted portions are not part of the corresponding generalized segments.}
	\label{fig:convergent}
\end{figure}

Two generalized segments in $\mathbb{R}^d$ are called \emph{convergent} if the two lines containing them meet at a point, that could be also the point at infinity of the two lines, that does not belong to any of the two generalized segments (see Figure \ref{fig:convergent}). In particular, two generalized segments on two lines that do not meet at any point are not convergent. 

\begin{theorem}
\label{t13}
Any set of $n$ red and blue points in $\mathbb{R}^d$ ($d \geq 2$) that is not contained
in a hyperplane determines at least $\frac{1}{48}dn-1$ generalized segments 
no two of which are convergent or collinear.
\end{theorem}

Notice that if we color all points of a set $P$ by red and because two parallel line segments are convergent as generalized segments, then 
Theorem \ref{t13} implies that a $d$-dimensional set of $n$ points in $\mathbb{R}^d$ determines at least $\frac{1}{48}dn-1$ lines with pairwise distinct directions.

When the $n$ points have the same color the notion of two segments being convergent coincides
with the same notion in \cite{PPS04} and \cite{PPS} which is that two segments are convergent if they are opposite edges in a convex quadrilateral. One of the main theorems in \cite{PPS04} is the following generalization of the above mentioned Ungar's theorem (\cite{Un}).

\begin{theorem}\label{theorem:pps}
A two dimensional set of $n$ points determines at least $2\lceil \frac{n}{2} \rceil$
segments no two of which are convergent.
\end{theorem}

We will use Theorem \ref{theorem:pps} to establish the basis of induction (on the dimension $d$) in the proof of Theorem \ref{t13}.

\vspace{0.2cm}

Using Theorem \ref{t12}, we obtain a new result about the distinct distances
problem for $d$-dimensional sets in typical $d$-norms. The
distinct distances problem in the Euclidean plane is one of the best known
classical open problems in Discrete Geometry, raised by Erd\H{o}s  in
1946 \cite{Er1}.  This is the problem of determining or estimating
the minimum possible number of distinct distances 
determined by $n$ points in the Euclidean plane. Although this
problem has been settled by Guth and Katz up to a 
$\sqrt {\log n}$ multiplicative factor \cite{GK}, the problem 
for $d$-dimensional Euclidean norms
for $d \geq 3$ remains wide open. Indeed, for each $d \geq 3$ the
conjecture is that the minimum possible number of distinct
distances determined by $n$ points in $\mathbb{R}^d$ with respect
to the Euclidean norm is $\Theta(n^{2/d})$. The $d$-dimensional 
integer box with edges of length $n^{1/d}$ shows that this is
an upper bound, but the best known lower bound does
not provide the correct exponent of $n$, see \cite{SV}.
The same problem for general norms has been considered as well,
see \cite{ABS} and the references therein. Call a norm $\|\cdot \|$ on
$\mathbb{R}^d$ a $d$-norm. A recent result proved in \cite{ABS}
asserts that there are $d$-norms in which any set
of $n>n_0(d)$ points determines at least $(1-o(1))n$ distinct
distances with respect to $\|\cdot\|$. 
In fact, this holds for all typical $d$-norms, where
the notion of a typical norm is defined as follows.

Identify a norm with its unit ball. The Hausdorff distance
between any two such unit balls $A$ and $B$ is the maximum
of all distances between a point of $A$ and the set $B$ and between a 
point of $B$ and the set $A$. This distance defines a metric, and
hence a topology, on the space of all $d$-norms. A set 
$S$ in this space is nowhere dense if every non-empty 
open set contains a nonempty open set which  does not intersect
$S$. A \emph{meagre} set is a countable union of nowhere dense sets. A
space is called a Baire space if the complement of each meagre set in it 
is dense. It is known that the space of all $d$-norms endowed with 
the Hausdorff metric as above is a Baire space. See e.g. \cite{ABS} 
for additional relevant
details and references. In this terminology it is proved in
\cite{ABS} that in all $d$-norms but a meagre set,
any set of $n>n_0(d)$ points determines at least $(1-o(1))n$
distinct distances, where the $o(1)$-term tends to $0$ as $n$ tends
to infinity. This is tight up to the $o(1)$-error term, 
as in every $d$-norm,
any set of $n$ points along an arithmetic progression on a line
determines exactly $n-1$ distinct distances. 

Note, however, that this extremal example that appears in any
$d$-norm is one-dimensional. It seems natural to consider the
distinct distances problem for configurations of $n$ points in
a $d$-norm that are $d$-dimensional, that is, do not all lie
in an affine hyperplane. For this case, we suggest the following
conjecture.

\begin{conjecture}
\label{c14}
For every fixed $d$ the following holds for all $d$-norms 
$\|\cdot \|$ but a
meagre set. For all $n>n_0(d)$, any set of $n$ points in 
$\mathbb{R}^d$ that
do not all lie in an affine hyperplane determine at least
$(d-o(1))n$ distinct distances with respect to $\|\cdot \|$, 
where the $o(1)$-term tends to $0$
as $n$ tends to infinity.
\end{conjecture}

This is, of course, trivial for $d=1$ but is open already for
$d=2$, as the result in \cite{ABS} only ensures $(1-o(1))n $
distinct distances for every fixed dimension $d$.

Here we observe that by combining the assertion of Theorem
\ref{t12} with the arguments in \cite{ABS} we get the following
weaker version of the conjecture.
\begin{theorem}
\label{t15}
There exists an absolute positive constant $b$ so that for
any fixed $d \geq 2$ and all
$d$-norms $\| \cdot \|$ but a meagre set the following holds.
Any set of $n>n_0(d)$ 
points in $\mathbb{R}^d$ that do not all lie in an affine hyperplane
determines at least $(bd-o(1))n$ distinct distances
with respect to $\|\cdot \|$, where 
the $o(1)$-term tends to $0$ as $n$ tends to infinity.
\end{theorem}

The rest of this note is organized as follows. In Section 2 we
describe the proof of Theorem \ref{t13}. Section 3 contains 
a sketch of the 
proof of Theorem \ref{t15}. The final section 4 contains some
concluding remarks and open problems.

\section{Distinct directions in $\mathbb{R}^d$}
In this section we prove Theorem \ref{t13}. We will need the following two dimensional lemma. We say that red is the opposite color of blue and vice versa. The idea behind the lemma and its proof are quite simple. Unfortunately, this simplicity is not quite reflected in the formulation of the lemma
that is a bit longer. We hope that the figures help to make the discussion a bit clearer.

\begin{lemma}\label{lemma:2d}
Let $x$ be a point that is either red or blue and let $P$
be a set of $n \geq 2$ red and blue points no two of which are collinear with $x$, where all points lie in the plane. In particular, the points of $P$ are all different from $x$. Let $\ell$ be a (horizontal) line through $x$ that does not pass through any point of $P$.
Partition $P$ into $P_{1} \cup P_{2}$ where $P_{1}$ is the set of points of $P$ above $\ell$ and $P_{2}$ is the set of points of $P$ below $\ell$. For every point $p$ in $P_{1}$ define $f(p)=p$.
In this case $f(p)$ has the same color as $p$. For every point $p \in P_{2}$ define $f(p)$ to be the reflextion of $p$ with respect to $x$. In this case we color $f(p)$ in the opposite color of $p$.
Then index the points of $P$ by $p_{1}, \ldots, p_{n}$ according to the (say increasing) angle between the segments $[x, f(p_{i})]$
and $\ell$ (see Figure \ref{fig:indexing}).

\begin{figure}[ht]
	\centering
	\includegraphics[width=6in]{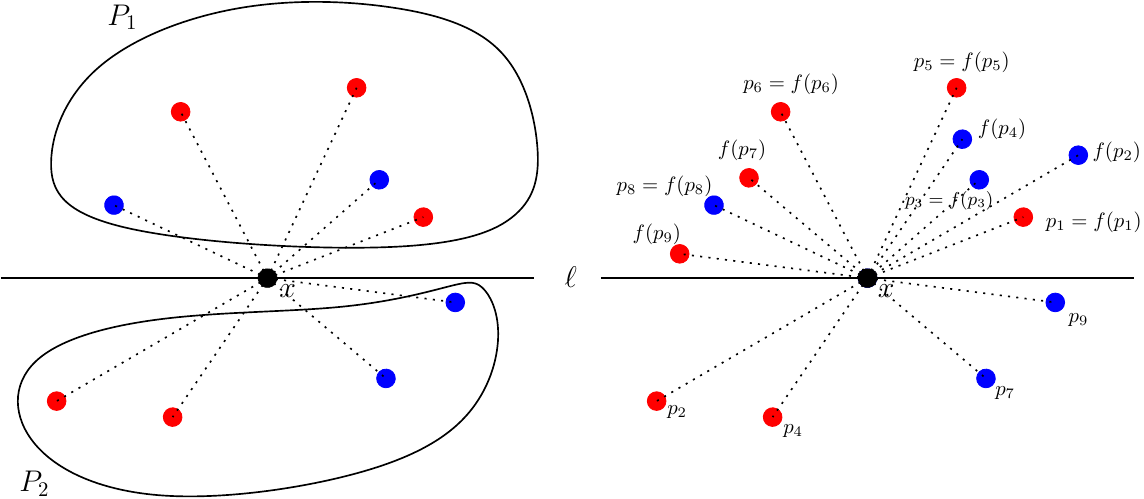}
	\caption[]
	{The indexing of the points of $P$ in Lemma \ref{lemma:2d}.}
	\label{fig:indexing}
\end{figure}

The following two assertions hold.

\noindent {\bf 1.} If $f(p_{1})$ and $f(p_{n})$ have the same color,
then the generalized segment determined by $p_{1}$ and $p_{n}$ is not convergent with respect to any of the generalized segments determined by $x$ and a point $p$ in $P$.

\medskip

\noindent {\bf 2.} For $1 \leq i < n$, if $f(p_{i})$ and $f(p_{i+1})$ have opposite colors, then the generalized segment determined by $p_{i}$ and $p_{i+1}$ is not convergent with respect to any of the generalized segments determined by $x$ and a point $p$ in $P$.
\end{lemma}

\noindent {\bf Proof.}
The proof is a fairly simple case analysis. We will try to shorten the analysis by using symmetry of cases.

\begin{figure}[ht]
	\centering
	\includegraphics[width=6in]{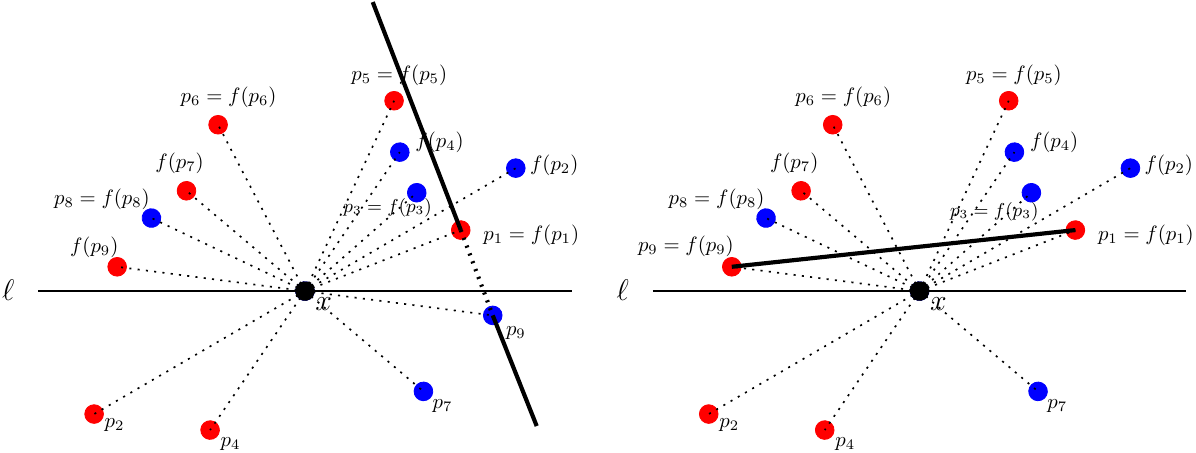}
	\caption[]
	{Part 1 of Lemma \ref{lemma:2d}. On the left the case where $p_{1}$ and $p_{n}$ (here $n=9$) are separated by $\ell$. On the right the case where $p_{1}$ and $p_{n}$ are not separated by $\ell$.}
	\label{fig:part1}
\end{figure}

We start by proving part 1 of the lemma (see Figure \ref{fig:part1}). If $p_{1}$ and $p_{n}$
are separated by $\ell$, then $p_{1}$ and $p_{n}$ have opposite colors and they appear cyclically consecutive around $x$
because of the way we indexed the points of $P$. The generalized segment determined by $p_{1}$ and $p_{n}$ is the line complement of the open interval $(p_{1},p_{n})$. Notice that no line through $x$ and a point of $P$ crosses this open interval. It follows now from the definition of convergent generalized segments that the generalized segment determined by $p_{1}$ and $p_{n}$ cannot be convergent with respect to any generalized segment determined by $x$ and a point of $P$.

If $p_{1}$ and $p_{n}$ are on the same side of $\ell$, then the generalized segment determined by $p_{1}$ and $p_{n}$ is the closed interval $[p_{1}, p_{n}]$. We notice that this interval is crossed by every line through $x$ and a point of $P$ because of the way we indexed the points of $P$. It follows now from the definition of convergent generalized segments that the generalized segment determined by $p_{1}$ and $p_{n}$ cannot be convergent with respect to any generalized segment determined by $x$ and a point of $P$.

\begin{figure}[ht]
	\centering
	\includegraphics[width=6in]{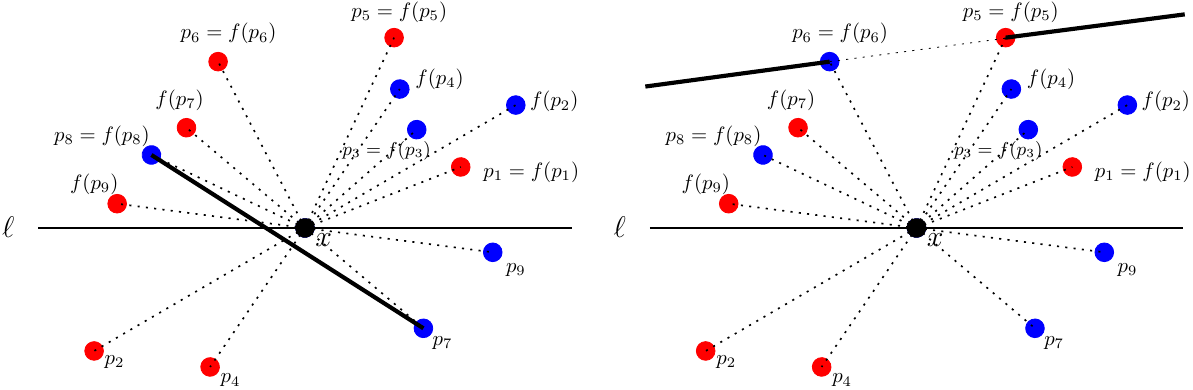}
	\caption[]
	{Part 2 of Lemma \ref{lemma:2d}. On the left the case where $p_{i}$ and $p_{i+1}$ (here $p_{7}$ and $p_{8}$) are separated by $\ell$. On the right the case where $p_{i}$ and $p_{i+1}$ (here $p_{5}$ and $p_{6}$) are not separated by $\ell$.}
	\label{fig:part2}
\end{figure}

Having proved part 1, we conclude by proving part 2 of the lemma (see Figure \ref{fig:part2}).
Assume first that $p_{i}$ and $p_{i+1}$ are not separated by $\ell$. Then by choosing $\ell$ to pass in such a way that it separates $p_{i}$ and $p_{i+1}$, we reduce the lemma to the case of part 1 that we proved already.

It is left to consider the case where $p_{i}$ and $p_{i+1}$ are separated by $\ell$. Then $p_{i}$ and $p_{i+1}$ have the same color. By choosing $\ell$ to pass in such a way that it separates $f(p_{i})$ and $f(p_{i+1})$, we reduce the lemma to the case of part 1 that we proved already.

\hfill $\Box$

\medskip

%{Proof of Theorem \ref{t13}}
Denote by $f(d,n)$ the maximum number such that any set of $n$ red and blue points that has affine dimension $d$ determines at least
$f(d,n)$ generalized segments no two of which are convergent.

We need to show that $f(d,n) \geq \frac{1}{48}{nd}-1$.
The proof goes by induction on $d$. The case $d=2$ follows from 
the so called weak Dirac's Theorem  (\cite{H17}) saying that in any noncollinear set $P$ of points in $\mathbb{R}^2$ one can find a point $x \in P$
that is incident to at least $\frac{n}{3}$ distinct lines through $x$ and another point of $P$. Indeed, notice that generalized segments sharing a common endpoint are by definition pairwise non-convergent.  In fact, this argument and the weaker lower bound of $\frac{n}{3}$ are valid in any dimension, showing that Theorem \ref{t13} is true for every $d<16$.
Assuming the statement is true for $d-1$, we prove it for $d$. 

We use a result in \cite{DSW14}, improving on a previous bound 
in \cite{BDWY11}. We say that a line is special with respect 
to a set of points $P$ if it contains at least 
$3$ points of $P$. We will use Theorem 1.9 in \cite{DSW14}: 

\begin{theorem}[Theorem 1.9 in \cite{DSW14}]
\label{theorem:DSW}
Let $P$ be a set of $n$ points with the property that for every 
$x \in P$ at least $\delta(n-1)$ of the rest of the points lie on 
special lines through $x$. Then the affine dimension of 
the set $P$ is at most $\frac{c}{\delta}$, where $c=12$.
\end{theorem}

An immediate consequence of Theorem \ref{theorem:DSW} is the following.

\begin{corollary}
\label{corollary:DSW}
Let $P$ be a $d$-dimensional set of $n$ points. Then there is a 
point $x$ in $P$ such that the number of lines connecting $x$ 
to the other points in $P$ is at least $n(1-\frac{c
}{d})$.
\end{corollary}

\noindent {\bf Proof.}
Denote by $M$ the maximum number such that there is a point in $P$
with $M$ distinct lines connecting it to the other point in $P$.
If $M=n-1$, there is nothing to prove. We therefore assume that $M<n-1$.

For every point $x \in P$ the number of points of $P \setminus \{x\}$ 
not lying on special lines through $x$ is at most $M-1$ (here we use 
the fact that $M<n-1$). Consequently, the number of points 
of $P \setminus \{x\}$ that do lie on special lines through $x$ is at least
$$n-1-(M-1)=n-M=\frac{n-M}{n-1}(n-1).$$ Taking $\delta=\frac{n-M}{n-1}$ 
and applying Theorem \ref{theorem:DSW}, we get 
$$
d \leq c \frac{n-1}{n-1-M}.
$$
This implies 
$$
M \geq (1-\frac{c}{d})(n-1)+1 \geq (1-\frac{c}{d})n.
$$
%\bbox
\hfill $\Box$

\medskip

We remark that the exact bound in Corollary \ref{corollary:DSW} 
is unknown and most likely difficult to achieve. See Section 4 for
more details.

Recall that $c=12$ in Theorem \ref{theorem:DSW}, however this value 
is not known to be tight and has been reported recently to be 
improved to $c=4$. For this reason we incorporate the constant 
$c$ in the rest of our calculations as a parameter rather than 
using its current best known published value $c=12$. We will also review the basis of induction.

We use Corollary \ref{corollary:DSW} to find a point $x \in P$ that we can connect by at least $n(1-\frac{c}{d})$ distinct lines to the other points in $P$. On each of these lines through $x$ we arbitrarily pick one of the points (if this line contains more than one point) different from $x$ to form a set $\tilde{P}$ of cardinality at least 
$n(1-\frac{c}{d})$. We define $S_{1}$ to be the collection of at least $n(1-\frac{c}{d})$ generalized segments determined by $x$ and a point of $\tilde{P}$. 

Next, we centrally project through $x$ the set $\tilde{P}$ onto a generic hyperplane $H$ of dimension $d-1$, not passing through $x$. On $H$ we get a set $P'$ of at least $n(1-\frac{c}{d})$ distinct points. Let $H'$ be the hyperplane parallel to $H$ that passes through $x$. We color a point $p'$ of $P'$ by the same color of its pre-image in $\tilde{P}$ if its pre-image is not separated by $H'$ from $H$. Otherwise we color $p'$ by the opposite color of its pre-image in $\tilde{P}$. 

\begin{figure}[ht]
	\centering
	\includegraphics[width=6in]{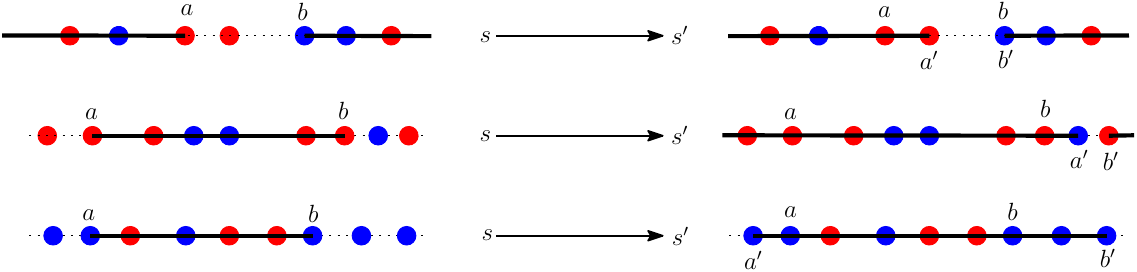}
	\caption[]
	{Replacement of the generalized segment $s$ by the generalized segment $s'$ that contains $s$. The black bold portions represent the respective generalized segments.}
	\label{fig:replacement}
\end{figure}

We notice that the set $P'$ is $(d-1)$-dimensional. This is because otherwise the set $P$ cannot be $d$-dimensional. 
By the induction hypothesis we find  $f(d-1, n(1-\frac{c}{d-1}))$ generalized segments
determined by points in $P'$ no two of which are convergent.
For each such generalized segment $s$ we replace $s$ by (possibly) another 
generalized segment $s'$ collinear with $s$ and that contains $s$ as follows (see Figure \ref{fig:replacement}). Let $\ell$ be the line that contains $s$. If $s$ is determined by two points $a$ and $b$ of $P'$ of opposite colors, we find two consecutive points $a'$ and $b'$ of $P'$ on $\ell$ that have different colors such that $a'$ and $b'$ are contained in $[a,b]$. We replace $s$ with the generalized segment $s'$ determined by $a'$ and $b'$. Notice that indeed $s$ was replaced by another
generalized segment $s'$ with endpoints of opposite colors and that $s'$ contains $s$. 

If $s$ is determined by two points $a$ and $b$ of $P'$ of the same color, we act differently. If there is a point of $P'$ on $\ell$ that is outside of $[a,b]$ and has opposite color of that of $a$ and $b$, we find on $\ell$ two consecutive points, $a'$ and $b'$, of $P'$ that lie outside of $(a,b)$ and have opposite colors. We then replace $s$ by the generalized segment $s'$ determined by $a'$ and $b'$. Notice that 
$s'$ contains $s$ because $s'$ is in fact the complement of $(a',b')$ on $\ell$. If, on the other hand, there is no point of $P'$ on $\ell$ outside of $[a,b]$ that has the opposite color of the color of $a$ and $b$, then we replace $s$ by the generalized segment $s'$ determined by the two extreme points of $P'$ on $\ell$. We notice that $s'$ contains $s$ also in this case.

We started with a collection of $f(d-1, n(1-\frac{c}{d-1}))$ generalized segments determined by points in $P'$ such that no two of the segments are convergent and we replaced every segment $s$ in this collection by another segment, collinear with $s$, that contains $s$. We therefore remain with a collection $S_{2}$ of $f(d-1, n(1-\frac{c}{d-1}))$ generalized segments determined by points in $P'$ such that no two of the segments are convergent.

Each segment $s \in S_{2}$ is determined by two points $a'$ and $b'$ of $P'$ that are the central projections through $x$ on $H$ of two points $\tilde{a}$ and $\tilde{b}$, respectively, in $\tilde{P}$.  We denote by $\tilde{s}$, the generalized segment determined by $\tilde{a}$ and $\tilde{b}$.
We set $\tilde{S_{2}}=\{\tilde{s} \mid s \in S_{2}\}$.

\begin{figure}[ht]
	\centering
	\includegraphics[width=6in]{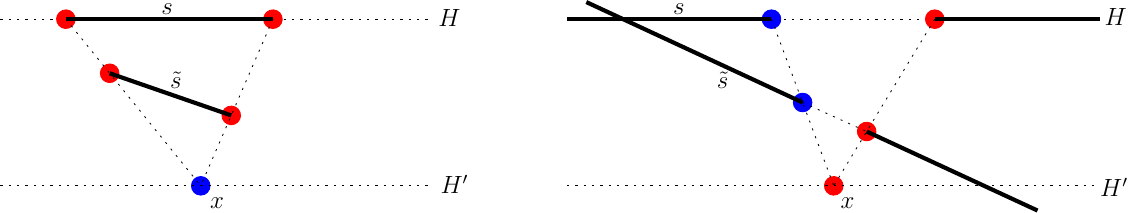}
	\caption[]
	{Claim \ref{claim:S12}. Illustrating the central projection through $x$ of $\tilde{s}$ on $H$ is equal to $s$ in the case the endpoints of $\tilde{s}$ are not separated by $H'$.}
	\label{fig:central_1}
\end{figure}

\begin{claim}\label{claim:S12}
$S_{1} \cup \tilde{S}_{2}$ is a collection of generalized segments no two of which are convergent and no two are collinear. 
\end{claim}

\noindent {Proof.}
Clearly, no two segments in $S_{1}$ are convergent because every  two contain the point $x$. We further notice that no two segments in $\tilde{S}_{2}$ are convergent. The crucial observation here is that every generalized segment $s \in S_{2}$ is equal to the central projection through $x$ of the generalized segment $\tilde{s}$. This is because of the way we color the points of $P'$. 

\begin{figure}[ht]
	\centering
	\includegraphics[width=6in]{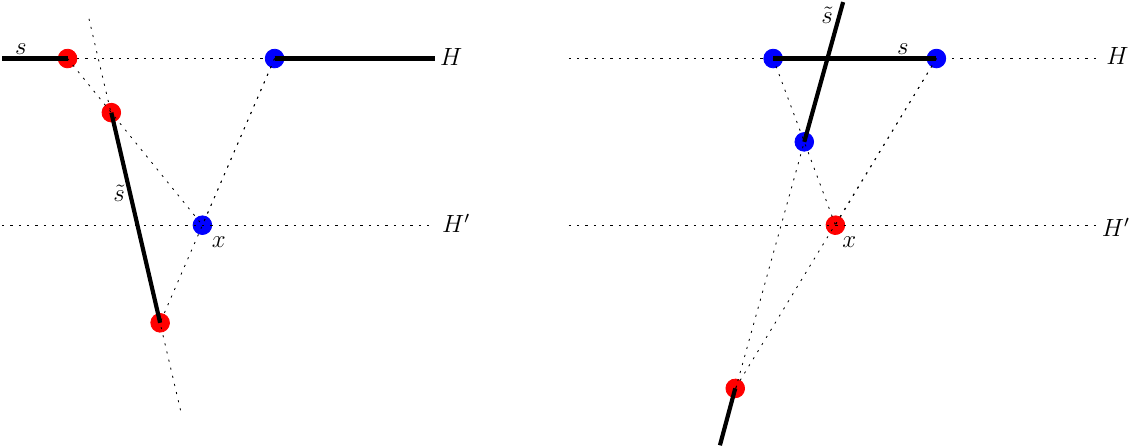}
	\caption[]
	{Claim \ref{claim:S12}. Illustrating the central projection through $x$ of $\tilde{s}$ on $H$ is equal to $s$ in the case the endpoints of $\tilde{s}$ are separated by $H'$.}
	\label{fig:central_2}
\end{figure}

One has to consider here the cases according to whether the 
two endpoints of $s$ have the same color or not and according to whether the two endpoints of $\tilde{s}$ are separated by $H'$ or not (see Figure \ref{fig:central_1} and Figure \ref{fig:central_2}). In each of these cases we come to the conclusion that the central projection, through $x$, of $\tilde{s}$ on $H$ is equal to $s$. 

The fact that no two generalized segments in $\tilde{S}_{2}$ are convergent follows now from the definition of convergent generalized segments and the fact that the central projection through $x$ is one-to-one on lines that do not pass through $x$.

It is left to show that it is not possible that a segment $s \in S_{1}$ and a segment $\tilde{s'} \in \tilde{S}_{2}$ are convergent.
This is precisely the content of Lemma \ref{lemma:2d}. 
More specifically, consider $s \in S_{1}$ and 
$\tilde{s'} \in \tilde{S}_{2}$. Assume to the contrary that they are convergent. Then there is a two dimensional plane $K$ containing both $s$ and $\tilde{s'}$. Let $\ell$ be the line of intersection of $K$ and $H'$ and apply Lemma \ref{lemma:2d}.
By our construction of the segments in $S_{2}$ it follows from Lemma \ref{lemma:2d} that $\tilde{s'}$ is non-convergent with any generalized segment on $K$ determined by $x$ and a point in $P'$.
In particular, $s$ and $\tilde{s'}$ are not convergent.

\hfill $\Box$

\medskip

Based on Claim \ref{claim:S12}, we can now write
$f(d,n) \geq n(1-\frac{c}{d})+f(d-1,n(1-\frac{c}{d}))$. 
We will now prove by induction on $d$ that $f(d,n) \geq \frac{dn}{4c}-1$. 
The constant $4$ can be improved, and we make no attempt to optimize
it here.

For the basis of induction, we use Theorem \ref{theorem:pps}. Given a set of $n$ red and blue points, we may assume without loss of generality that $\frac{n}{2}$ of them are of the same color red. If this set of red points is not collinear, then Theorem \ref{theorem:pps} implies, after a generic projection to a $2$-dimensional plane, that the red points determine at least $\frac{n}{2}-1$ (generalized) segments, no two of which are convergent. Therefore, $f(d,n) \geq \frac{n}{2}-1$. We come to the same conclusion
if the set of red points is collinear by choosing a (blue) point $x$ not collinear with the red points and then the generalized segments determined by $x$ and one of the red points is a collection of at least $\frac{n}{2}$ segments, no two of which are convergent.

It follows now that the inequality $f(d,n) \geq \frac{dn}{4c}-1$ is true for $d \leq 2c$ 
because, as we have seen, $f(d,n) \geq \frac{n}{2}-1$ for any $d \geq 2$.

Assume therefore that $d>2c$. We prove 
by induction on $d$ that $f(d,n) \geq \frac{dn}{4c}-1$. 
For the induction step, we need the following inequality to hold.

$$
n(1-\frac{c}{d})+\frac{(d-1)
n(1-\frac{c}{d})}{4c}-1 \geq \frac{dn}{4c}-1.
$$

After simplification, we get the following equivalent inequality to prove:

$$
\frac{3}{4}+\frac{1}{4d} \geq \frac{1}{4c}+\frac{c}{d}
$$

As we assume $d \geq 2c$, it is enough to show

$$
\frac{1}{4} +\frac{1}{4d}\geq \frac{1}{4c}
$$

One can easily check that this inequality holds for all $c \geq 1$. 
In our case $c=12$ and we have therefore proved that 
$f(d,n) \geq \frac{dn}{48}-1$. 
This concludes the proof of Theorem \ref{t13}, which implies
Theorem \ref{t12}.

%The inequality $f(d,n) \geq \frac{dn}{4c}$ is true for $d \leq 3c$ 
%because of Ungar's theorem, that is $f(d,n) \geq f(2,n) \geq n-1$.

%Assume therefore that $d>3c$. We prove 
%by induction on $d$ that $f(d,n) \geq \frac{dn}{4c}$. 
%For the induction step we need the following inequality to hold.

%$$
%n(1-\frac{c}{d})+\frac{(d-1)
%n(1-\frac{c}{d})}{4c} \geq \frac{dn}{4c}.
%$$

%After simplification we get the following equivalent inequality to prove:

%$$
%\frac{3}{4}+\frac{1}{4d} \geq \frac{1}{4c}+\frac{c}{d}
%$$

%As we assume $d \geq 3c$, it is enough to show

%$$
%\frac{5}{12} +\frac{1}{4d}\geq \frac{1}{4c}
%$$

%One can easily check that this inequality holds for all $c \geq 1$. 
%In our case $c=12$ and we have therefore proved that 
%$f(d,n) \geq \frac{dn}{48}$. 
%This concludes the proof of Theorem \ref{t13}, which implies
%Theorem \ref{t12}.

\section{Distinct distances for typical $d$-norms}
In this section we sketch the proof of Theorem \ref{t15}. 
The proof follows
closely the one in \cite{ABS}, replacing Ungar's Theorem stated 
as Theorem 4.6 in \cite{ABS} by our Theorem \ref{t12} here. Since
the proof is very similar to the one in \cite{ABS} we do not repeat
here the identical parts, and merely describe the different points,
referring to the arguments in \cite{ABS} whenever needed.
We start with the following  modified version of Lemma 4.2 in
\cite{ABS}.

\begin{lemma}
\label{l31}
There exists an absolute positive constant $b$ so that the
following holds.
Let $d\ge 1$ be an integer and let $0<\mu<1$. 
Suppose that $n$ is sufficiently large with respect to $d$ and $\mu$. 
Let $F\su \mathbb{R}$ be a subfield of $\mathbb{R}$, 
and let $V$ be a vector space over $\mathbb{R}$. 
Let $\mathbf{u}_1,\dots, \mathbf{u}_k\in V$ be non-zero 
vectors in $V$, and let $\mathbf{p}_1,\dots,\mathbf{p}_n\in V$ 
be distinct vectors such that not all 
of $\mathbf{p}_1,\dots,\mathbf{p}_n$ lie in a common
$(d-1)$-dimensional affine subspace of $V$ 
(as a vector space over $\mathbb{R}$). 
Suppose that for all $x,y\in\{1,\dots,n\}$ 
we have $\mathbf{p}_x-\mathbf{p}_y\in \mbox{span}_F(\mathbf{u}_i)$ 
for some $i\in \{1,\dots,k\}$. Then there exists a subset 
$I\su \{1,\dots,k \}$, such that we have 
$\mathbf{u}_\ell\in \mbox{span}_{F}(\mathbf{u}_i\midd i\in I)$ 
for at least $d\cdot |I|+(bd-\mu)\cdot n+1$ 
indices $\ell\in \{1,\dots,k     \}$.
\end{lemma}

\noindent
{\bf Proof (sketch).}
Setting $m=\lceil (bd-\mu)\cdot n\rceil\le bdn$, where $b$ is the
constant from Theorem \ref{t12},
we want to prove that there is a subset $I\su \{1,\dots,k\}$ 
with $\mathbf{u}_\ell\in \mbox{span}_{F}(\mathbf{u}_i\midd i\in I)$ 
for at least $d\cdot |I|+m+1$ indices $\ell\in \{1,\dots,k\}$. 
Suppose towards a contradiction that the desired subset 
$I\su \{1,\dots,k\}$ does not exist. Then for every subset 
$I\su \{1,\dots,k\}$, we have 
$\mathbf{u}_\ell\in \mbox{span}_{F}(\mathbf{u}_i\midd i\in I)$ 
for at most $d\cdot |I|+m$ indices $\ell\in \{1,\dots,k\}$.

We may assume that for every $i\in \{1,\dots,k\}$ there exist 
distinct $x,y\in \{1,\dots,n\}$ with $\mathbf{p}_x-\mathbf{p}_y\in 
\mbox{span}_{F}(\mathbf{u}_i)$ since otherwise, we can just omit 
all indices $i$ for which this is not the case, and relabel 
the remaining indices.

By Theorem \ref{t12}, there are at least $bdn$ different line 
directions in $\mbox{span}_\R(\mathbf{p}_1,\dots,\mathbf{p}_n)\su V$ 
appearing among the differences $\mathbf{p}_x-\mathbf{p}_y$ 
with $1\le x<y\le n$. For each of these differences we have 
$\mathbf{p}_x-\mathbf{p}_y\in \mbox{span}_F(\mathbf{u}_i)$ 
for some $i\in \{1,\dots,k\}$.
Hence, there must be at least $bdn$ 
different vectors $\mathbf{u}_i$, so $k\ge bdn$.

We now construct a sequence of distinct indices 
$j_1,\dots,j_r\in \{1,\dots,k\}$ recursively exactly as described
in the proof of Lemma 4.2 in \cite{ABS}. We also define the subsets
$H_0, H_1, \ldots ,H_r$ as in this proof and observe that the
following two claims hold just as in that proof:

\begin{claim}
\label{cl32} We have $r\le \frac{\mu}{3d}\cdot n$.
\end{claim}

\begin{claim}
\label{cl33} We have $|H_r|> \frac{\mu}{3d}\cdot n$.
\end{claim}

The only required difference between 
the proof of Claim \ref{cl33} here and that of Claim 2 in 
\cite{ABS} is the replacement of the
penultimate line of this proof which is

\[k\ge n-1\ge \frac{2\mu}{3}\cdot n+\lceil (1-\mu)\cdot 
n\rceil+1 = d\cdot \frac{2\mu}{3d}\cdot n+m+1\ge d\cdot |I|+m+1,\]

by the line

\[k\ge bdn \ge \frac{2\mu}{3}\cdot n+\lceil (bd-\mu)\cdot 
n\rceil+1 = d\cdot \frac{2\mu}{3d}\cdot n+m+1\ge d\cdot |I|+m+1,\]

\noindent
which holds 
if $n$ is sufficiently large with respect to $\mu$. 
This contradicts our assumption that such a set $I$ does not exist
and completes the proof of the claim.

The rest of the proof of the lemma is identical to that of Lemma
4.2 in \cite{ABS} where the only differences are as follows.
\begin{itemize}
\item
In the inequality

\[\sum_{j\in J}\lambda_j\ge \frac{\lambda_1+\dots+\lambda_k-
m\cdot \frac{3d}{\mu}}{d}\ge \frac{\frac{\mu}{24d}\cdot 
|H_r|^2-n\cdot \frac{3d}{\mu}}{d}\]

the last $n$ in the numerator in the right-hand-side has to be
replaced by $bdn$ giving

\[\sum_{j\in J}\lambda_j\ge \frac{\lambda_1+\dots+\lambda_k-
m\cdot \frac{3d}{\mu}}{d}\ge \frac{\frac{\mu}{24d}\cdot 
|H_r|^2-bdn\cdot \frac{3d}{\mu}}{d}.\]
\item
In the next sentence $n$ has to be  replaced, again, by $bdn$ as
written in the following three lines:

\noindent
By Claim \ref{cl33} 
we have $bdn\cdot \frac{3d}{\mu}\le bd 
\left(\frac{3d}{\mu}\right)^2\cdot |H_r|\le \frac{\mu}{48d}
\cdot |H_r|^2$ if $n$ (and therefore also $|H_r|\ge 
\frac{\mu}{3d}\cdot n$) is sufficiently large 
with respect to $d$ and $\mu$. So we can conclude that

\[\sum_{j\in J}\lambda_j\ge\frac{\frac{\mu}{24d}\cdot |H_r|^2
-bd n\cdot \frac{3d}{\mu}}{d}\ge \frac{\frac{\mu}{48d}
\cdot |H_r|^2}{d}=\frac{\mu}{48d^2}\cdot |H_r|^2.\]
\end{itemize}

The end of the proof of the lemma is identical to the one in \cite{ABS}. 
%\bbox
\hfill $\Box$
\vspace{0.2cm}

\noindent
The proof of Theorem \ref{t15} is now essentially identical to the
proof of Theorem 1.3 as described in Section 5 of \cite{ABS}, where
the only differences are as follows.
\begin{itemize}
\item
Each of the appearances of the
quantity $(1-\mu)n$ should be replaced by $(bd-\mu))n$, and in
particular the 
parameter $m$ which is $\lceil (1-\mu)n \rceil$ in Section 5
of \cite{ABS} has to
be $\lceil (bd-\mu)n \rceil$ here.
\item
The paragraph considering the case that all the points 
$\bf p_1, \ldots ,\bf p_n$ are on a common line should be omitted
here,
as in our case, by assumption, the affine dimension of this set of
points is $d \geq 2$.
\end{itemize}
This completes the sketch of the proof of Theorem \ref{t15}.
%\bbox
\hfill $\Box$

\section{Concluding remarks and open problems}

\begin{itemize}
\item
The discussion in Section 3 shows that while
Conjecture \ref{c11}, if true, does not imply Conjecture \ref{c14}
as stated here, it does imply that any $d$-dimensional set of
$n>n_0(d)$ in a typical $d$-norm determines at least
$(d-1-o(1))n$ distinct distances.

The argument together with the main result of \cite{PPS} also shows
that any $3$-dimensional set of $n>n_0$ points in a typical $3$-norm
determines at least $(2-o(1))n$ distinct distances.
\item
A slightly modified version of an old conjecture of Dirac \cite{Di}
asserts that for any set of $n$ points in the plane, not all on a 
single line, there is a point that lies in at least 
$n/2-O(1)$ distinct lines determined by the set. This Conjecture
is still open but several weaker versions have been established
over the years. See \cite{H17} for the best known bound and the
history of the problem.
Corollary \ref{corollary:DSW} deals with the higher dimensional
version of the problem. Its quality depends on the best 
known estimate for the constant $c$. It is worth noting that
up to the constant $c$ the bound it provides is 
essentially best possible. Indeed, 
to see an example
assume for simplicity that $d$ is odd and $n$ is divisible 
by $d+1$. Then consider $(d+1)/2$ skew lines
that affinely span $\mathbb{R}^d$. On each of these lines 
take $2n/(d+1)$ points. It is easy to check that for every point 
of the resulting
set of points, the number of lines connecting it 
to the other points of the set is $n-2n/(d+1)+1=n(1-\frac{2}{d+1})+1$.
\item
Theorem \ref{t15} together with the simple example of
$n-d+1$ points in an arithmetic progression on a line and $d-1$
additional points in general position not on this line, 
show that for typical $d$-norms, the minimum
possible number of distinct distances determined by a
$d$-dimensional set of $n>n_0(d)$ points is $\Theta(nd)$. Note that
this expression is an increasing function of $d$.
This is in contrast to the behavior of this minimum 
in Euclidean spaces, where the known upper bound 
(which is conjectured to be tight)
is $O(n^{2/d})$-a decreasing function of $d$.
\end{itemize}

\bigskip

\noindent {\bf \Large Acknowledgements.}
We thank Boyan Duan, Minghui Ouyang, and Zheng Wang for pointing out a gap in an earlier version of this paper. Because of this gap that they pointed out we had to use the definition of generalized segments and prove Lemma \ref{lemma:2d}, thus getting the stronger bichromatic version of Theorem \ref{t13}. Many Thanks !

\end{document}